\documentclass[aos,preprint]{imsart}
\RequirePackage{amsmath,amsfonts}
\RequirePackage[authoryear]{natbib}
\RequirePackage[colorlinks,citecolor=blue,urlcolor=blue,linkcolor=blue]{hyperref}

\startlocaldefs
\numberwithin{equation}{section}
\newcommand{\R}{\mathbb{R}}
\newcommand{\C}{\mathbb{C}}
\newcommand{\E}{\mathbb{E}}
\newcommand{\Pp}{\mathbb P}
\newcommand{\norm}[1]{\left\lVert#1\right\rVert}

\newtheorem{theorem}{Theorem}[section]
\newtheorem{lemma}[theorem]{Lemma}
\newtheorem{proposition}[theorem]{Proposition}
\newtheorem{corollary}[theorem]{Corollary}
\endlocaldefs

\begin{document}
\begin{frontmatter}
\title{Sharp Hellinger Risk of the Multivariate Gaussian NPMLE}
\runtitle{Sharp Hellinger Risk of the Gaussian NPMLE}
\runauthor{Wang and Li}

\begin{aug}
\author[A]{\fnms{Xingzhou}~\snm{Wang}%
  \ead[label=e1]{xwang268@ncsu.edu}}
\author[B]{\fnms{Zheng}~\snm{Li}%
  \ead[label=e2]{zli42@ncsu.edu}}

\address[A]{Department of Statistics, North Carolina State University%
  \printead[presep={,\ }]{e1}}
\address[B]{Department of Agricultural and Resource Economics, North Carolina State University%
  \printead[presep={,\ }]{e2}}
\end{aug}

\begin{abstract}
For multivariate Gaussian location mixtures whose mixing distributions
have coordinatewise sub-Gaussian tails, we prove that the unconstrained
nonparametric maximum likelihood estimator attains the minimax rate
$(\log n)^d/n$ under squared Hellinger loss. This holds for every fixed
dimension $d\geq2$ and every fixed tail parameter $\sigma>0$.
The upper bound applies to every measurable choice of a global
likelihood maximizer and improves the earlier bound of order
$(\log n)^{d+1}/n$ by a factor of $\log n$.

A key step is a sharp bracketing bound for multivariate polynomials.
For the $L^2(\rho)$ unit ball of real polynomials in $d$ variables
with total degree at most $m$, the log bracketing number at width
$\varepsilon$ is at most $\binom{m+d}{d}\log(C_d/\varepsilon)$
for $0<\varepsilon\leq1$. The constant $C_d$ depends only on $d$.
The bound holds for every positive Borel measure $\rho$ under which
each nonzero polynomial in the space has a finite, positive
$L^2(\rho)$ norm. We prove the bound using volume estimates for a
complex Veronese cone and Gaussian small-ball probabilities, then
apply it to bound the local Hellinger bracketing entropy of Gaussian
mixture densities. We obtain the matching minimax lower bound by adapting the Hermite perturbation construction of Kim and Guntuboyina.
\end{abstract}

\begin{keyword}[class=MSC]
\kwdgroup[type=primary]{\kwd{62G07}}
\kwdgroup[type=secondary]{\kwd{62G05}}
\end{keyword}
\begin{keyword}
\kwd{Gaussian location mixtures}
\kwd{nonparametric maximum likelihood}
\kwd{minimax risk}
\kwd{Hellinger distance}
\kwd{bracketing entropy}
\end{keyword}
\end{frontmatter}

\section{Introduction}\label{sec:intro}
Let $\phi_d$ be the $N(0,I_d)$ density. Write $\mathcal P(S)$ for the
set of Borel probability measures supported on $S$. We observe independent
random vectors $X_1,\ldots,X_n$ with common density
\begin{equation}\label{eq:mix-model}
 f_G(x) = \int_{\R^d}\phi_d(x-\theta)\,dG(\theta),
 \qquad G\in\mathcal P(\R^d).
\end{equation}
We estimate the mixture density $f_G$ under the squared Hellinger loss
$H^2(f,g)=\int_{\R^d}(\sqrt f-\sqrt g)^2\,dx$.

The unconstrained nonparametric maximum likelihood estimator (NPMLE)
is any measurable map $\widehat G_n$ satisfying
\begin{equation}\label{eq:mle}
 \widehat G_n\in
 \operatorname*{arg\,max}_{Q\in\mathcal P(\R^d)}
 \sum_{i=1}^n\log f_Q(X_i).
\end{equation}
In particular, the maximization in \eqref{eq:mle} uses neither a known
support bound nor the tail parameter of $G$.
Lemma~\ref{lem:convex-hull-support} shows that a maximizer exists for
every sample and that we can select one measurably.

For $\sigma>0$, let $\mathcal S_d(\sigma)$ be the class of mixing laws
$G$ that have a mean $a\in\R^d$ and satisfy
\begin{equation}\label{eq:tailclass}
 G\{|\theta_j-a_j|>t\}\leq 2\exp\{-t^2/(2\sigma^2)\},
 \qquad t>0,\quad j=1,\ldots,d .
\end{equation}
The mean is unrestricted. We hold $d$ and $\sigma$ fixed as $n$ grows.
Write $\log_+ t=\max\{0,\log t\}$.

For sub-Gaussian mixing distributions, \citet[Corollary~2.2]{SG2020}
give an upper bound of order $(\log n)^{d+1}/n$ for the squared
Hellinger risk of the unconstrained NPMLE.
\citet[Theorem~2.2]{KG2022} give a minimax lower bound of order
$(\log n)^d/n$ on a sub-Gaussian class.

In one dimension, \citet[Theorem~1.1(i)]{HM2026} show that the
unconstrained NPMLE attains the minimax rate $(\log n)/n$ for each fixed
$\sigma\geq1$. In Section~5, they ask whether the unconstrained NPMLE
also attains the minimax rate for $d\geq2$. They also conjecture
a sharp order for the $L^2$ bracketing entropy of multivariate
polynomial balls. We answer both questions.

\begin{theorem}\label{thm:main}
Fix an integer $d\geq1$ and $\sigma>0$. For every $n\geq2$, let
$\widehat G_n$ be any measurable NPMLE satisfying \eqref{eq:mle}.
Then
\begin{equation}\label{eq:main}
 \sup_{G\in\mathcal S_d(\sigma)}
 \E_{f_G^{\otimes n}}H^2(f_{\widehat G_n},f_G)
 \leq C_{d,\sigma}\min\!\left\{1,\frac{(\log n)^d}{n}\right\},
\end{equation}
where $0<C_{d,\sigma}<\infty$ depends only on $d$ and $\sigma$.
\end{theorem}

\begin{corollary}\label{cor:minimax}
Fix an integer $d\geq2$ and $\sigma>0$. Let $\widehat G_n$ be any
measurable NPMLE satisfying \eqref{eq:mle}. There are constants
$0<c_{d,\sigma}\leq C_{d,\sigma}<\infty$ such that, for all sufficiently
large $n$,
\begin{align}
 c_{d,\sigma}\frac{(\log n)^d}{n}
 &\leq \inf_{\widetilde f_n}
       \sup_{G\in\mathcal S_d(\sigma)}
       \E_{f_G^{\otimes n}}H^2(\widetilde f_n,f_G) \notag\\
 &\leq \sup_{G\in\mathcal S_d(\sigma)}
       \E_{f_G^{\otimes n}}H^2(f_{\widehat G_n},f_G) \notag\\
 &\leq C_{d,\sigma}\frac{(\log n)^d}{n}.
 \label{eq:minimax}
\end{align}
The infimum is over all measurable density estimators. If $\sigma\geq1$,
we can choose the lower-bound constant to depend only on $d$.

The same bounds hold if we replace \eqref{eq:tailclass} by the condition
\[
 G\bigl\{|\langle u,\theta-a\rangle|>t\bigr\}
 \leq 2\exp\!\left(-\frac{t^2}{2\sigma^2}\right)
 \qquad\text{for all }t>0\text{ and all unit vectors }u\in\R^d.
\]
\end{corollary}

The matching lower bound uses the Hermite perturbation construction of
\citet[Section~3.3]{KG2022}. Appendix~\ref{sec:lower} adjusts the perturbation scales and normalization. It verifies that the resulting mixing laws satisfy the directional tail condition in Corollary~\ref{cor:minimax} for every fixed $\sigma>0$.

\subsection{Proof strategy and relation to earlier results}
The space of polynomials in $d$ variables of total degree at most $m$
has dimension $N=\binom{m+d}{d}$. Write $N_{[\,]}(\varepsilon)$ for the
$L^2$ bracketing number of its unit ball. In an $L^2$-orthonormal basis, a Euclidean $\delta$-net of the
coefficient unit ball gives brackets of $L^2$ width $2\delta\sqrt N$.
We can choose the net with at most $(1+2/\delta)^N$ points.
Taking $\delta=\varepsilon/(2\sqrt N)$ gives
\[
  \log N_{[\,]}(\varepsilon)
  \leq N\log\!\left(1+\frac{4\sqrt N}{\varepsilon}\right).
\]

Following \citet[Section~4]{HM2026}, we use a norm based on normalized
evaluation vectors and compare volumes in coefficient space. In one
dimension, these vectors trace a curve. Hundrieser and Manole use its
length to obtain a lower bound on the volume of the unit ball of this norm.

For multivariate polynomials, we show that the normalized evaluation
vectors lie in a linear image of the complex Veronese cone. To bound the covering number of the evaluation set, we compare the total volume of the cone inside disjoint small Euclidean balls with its volume inside one larger ball. We then use a Gaussian small-ball
estimate to obtain a lower bound on the volume of the unit ball of
the evaluation norm.
Theorem~\ref{thm:poly} gives
\[
 \log N_{[\,]}(\varepsilon)
 \leq N\log\!\left(\frac{C_d}{\varepsilon}\right),
 \qquad 0<\varepsilon\leq1,
\]
where $C_d$ depends only on $d$.

For mixing laws supported on $[-R,R]^d$, we adapt the weighted
decomposition of \citet[Section~3.2]{HM2026} to express the
difference between the square roots of two mixture densities.
We divide this difference by the radius $s$ of the local Hellinger
ball. On a larger cube in the observation space, we approximate
the factor multiplying the weight by a polynomial of total degree
$O_d(1+R^2+\log_+(1/s))$.

We use a grid on this cube to cover the weights by finitely many
functions at a fixed multiplicative accuracy. In one dimension, \citet[Proposition~3.8]{HM2026} give a similar cover for finite Gaussian mixtures. We combine our weight cover with the polynomial bracketing bound in
Theorem~\ref{thm:poly} to construct brackets of width at most
$s/2$. The same bound on the number of brackets holds for
every center of the local Hellinger ball.

We iterate this local bracketing construction at radii
$s,s/2,s/4,\ldots$. We use the resulting entropy bounds to control the local bracketing entropy integral. We then apply a likelihood risk bound. For $R\geq1$, the NPMLE over $\mathcal P([-R,R]^d)$ has squared
Hellinger risk at most
\[
 C_d\min\!\left\{1,\frac{(1+R^2+\log n)^d}{n}\right\},
\]
uniformly over all true mixing laws in $\mathcal P([-R,R]^d)$.
The constant $C_d$ depends only on $d$. We truncate the true mixing law and use the fact that every
unconstrained maximizer is supported on the convex hull of the sample. We take $R$ of order $\sqrt{\log n}$ for fixed $d$ and $\sigma$.
This gives the claimed risk bound for the unconstrained NPMLE defined in
\eqref{eq:mle}. In one dimension, \citet[Remark~2.8]{HM2026} note that one can also derive the risk bound for the unconstrained NPMLE by approximating
the true sub-Gaussian mixing law with one supported on $[-R,R]$.

For univariate Gaussian location mixtures with sub-Gaussian mixing
laws, Theorem~1 of \citet{Zhang2009} implies an $O((\log n)^2/n)$
bound on the squared Hellinger risk of the unconstrained NPMLE. In dimension $d$, \citet[Corollary~2.2, part~2]{SG2020} give the corresponding bound $O((\log n)^{d+1}/n)$.
\citet{SoloffEtAl2025} study the NPMLE for multivariate Gaussian location mixtures with known covariance matrices that may vary across observations.
For Gaussian mixtures with mixing laws supported on a fixed compact
set, \citet[Corollary~11]{JiaEtAl2023} characterize the minimax squared
Hellinger risk in terms of local Hellinger covering numbers.

When the true mixing law has compact support, \citet[Theorem~2.1(i)]{GhoshEtAl2026} give a high-probability
bound on the squared Hellinger error of approximate maximizers of the unconstrained likelihood. When the true and fitted mixing laws are supported on the same
fixed, known compact set, \citet[Theorem~1(iv) and
Proposition~4]{Jiang2026} gives a high-probability upper bound
on the squared Hellinger error of order
\[
 \frac{(\log n)^{d+1}}{n(\log\log n)^d}
\]
for the ordinary constrained NPMLE and a suitably reweighted NPMLE. The reweighted estimator uses i.i.d.\ Gamma weights
drawn independently of the data, with the concentration parameter specified in equation~(6) of that paper.
\citet[Theorem~1.1]{HanEtAl2026} bound the number of atoms in the NPMLE when the true and fitted mixing laws are supported on the same fixed cube. \citet{MaEtAl2025} give bounds on the number of components needed to approximate Gaussian location mixture densities to a given accuracy. In one dimension, \citet[Theorem~19]{ChenWu2026} give an $O((\log n)/n)$ bound on the squared Hellinger risk of the support-constrained NPMLE.

Our risk bound holds for every measurable global maximizer of the log likelihood in \eqref{eq:mle}. \citet{YWR2024} propose a Wasserstein--Fisher--Rao gradient descent algorithm for this maximization problem. \citet[Theorem~2.1]{WFR2026} prove that the full sequence of explicit split iterates converges weakly throughout the original step-size range.

\section{Polynomial bracketing via complex geometry}
\label{sec:poly}

Let $m\geq 0$ be an integer. Write $\mathcal P_m(\R^d)$ for the space
of real polynomials in $d$ variables of total degree at most $m$.
Its dimension is
\[
  N=\dim\mathcal P_m(\R^d)=\binom{m+d}{d}.
\]

Let $l$ and $u$ be measurable real-valued functions on $\R^d$ with
$l\leq u$. Define the bracket
\[
  [l,u]
  =\{f\text{ measurable}: l(x)\leq f(x)\leq u(x)
       \text{ for every }x\in\R^d\}.
\]
Its $L^2(\rho)$ width is $\norm{u-l}_{L^2(\rho)}$.
The bracketing number $N_{[\,]}(\varepsilon,\mathcal F,L^2(\rho))$
is the smallest number of brackets of $L^2(\rho)$ width at most
$\varepsilon$ needed to cover $\mathcal F$.

\begin{theorem}[Polynomial bracketing]
\label{thm:poly}
Fix integers $d\geq1$ and $m\geq0$. Let $\rho$ be a positive Borel
measure on $\R^d$ with
\[
  \int_{\R^d}(1+\norm{x}_2^{2m})\,d\rho(x)<\infty.
\]
Assume that $\norm{p}_{L^2(\rho)}>0$ for every nonzero
$p\in\mathcal P_m(\R^d)$. Set
\[
  N=\binom{m+d}{d},
  \qquad
  \mathcal F_m(\rho)=
  \{p\in\mathcal P_m(\R^d):\norm{p}_{L^2(\rho)}\leq1\}.
\]
There is a constant $C_d\geq1$ depending only on $d$ such that, for
$0<\varepsilon\leq1$,
\begin{equation}
\label{eq:poly-entropy}
  \log N_{[\,]}(\varepsilon,\mathcal F_m(\rho),L^2(\rho))
  \leq N\log\!\left(\frac{C_d}{\varepsilon}\right).
\end{equation}
For $0<\varepsilon\leq1$, we also have
\[
  \log N_{[\,]}(\varepsilon,\mathcal F_m(\rho),L^2(\rho))
  \geq N\log(1/\varepsilon).
\]
\end{theorem}

If $\rho$ has a Lebesgue density that is positive almost everywhere
on a nonempty open cube $U$, then $\norm{p}_{L^2(\rho)}>0$ for
every nonzero $p\in\mathcal P_m(\R^d)$.
If $\norm{p}_{L^2(\rho)}=0$, then $p=0$ almost everywhere on $U$.
By continuity, $p(x)=0$ for every $x\in U$.
A polynomial that vanishes on a nonempty open set is identically zero.

Following \citet[Section~4]{HM2026}, we use an evaluation norm
to compare volumes in coefficient space. In the multivariate case, we give a lower bound on the volume of
the unit ball of this norm, uniformly in $\rho$. We first bound
the Euclidean covering numbers of the set of normalized evaluation
vectors. We compare the total volume of the complex cone inside disjoint
small balls with its volume inside a larger ball. This gives a
covering bound with a factor of $m^d$. For comparison, the general
covering bounds for real algebraic sets and polynomial images
summarized in \citet[Table~1]{ZK2026} depend explicitly on the ambient dimension.

We assume $m\geq1$ for now and treat $m=0$ separately.
Choose an $L^2(\rho)$-orthonormal basis $q_1,\ldots,q_N$ of
$\mathcal P_m(\R^d)$. For $x\in\R^d$ and $a\in\R^N$, define
\[
  Q(x)=(q_1(x),\ldots,q_N(x)),\qquad
  E(x)=\norm{Q(x)}_2,\qquad
  p_a(x)=\langle a,Q(x)\rangle.
\]
The constant polynomial $1$ lies in the span of $q_1,\ldots,q_N$.
So $Q(x)\neq0$ and $E(x)>0$ for every $x\in\R^d$. Define
\begin{equation}
\label{eq:poly-V-star}
\begin{split}
  V&=\overline{\{Q(x)/E(x):x\in\R^d\}}\subset S^{N-1},\\
  \norm{a}_*
    &=\sup_{v\in V}|\langle a,v\rangle|
     =\sup_{x\in\R^d}\frac{|p_a(x)|}{E(x)},\\
  B_*&=\{a\in\R^N:\norm{a}_*\leq1\}.
\end{split}
\end{equation}
The set $V$ is compact because it is closed in $S^{N-1}$.
We use Euclidean distance on spheres throughout and write $B(v,r)$
for the open Euclidean ball of radius $r$ centered at $v$.

\begin{lemma}[Evaluation norm]
\label{lem:poly-evaluation}
The function $\norm{\cdot}_*$ is a norm on $\R^N$. Let $B_2^N$ denote the closed Euclidean unit ball. Then
\begin{equation}
\label{eq:poly-ball-inclusions}
  B_2^N\subseteq B_*\subseteq \sqrt N\,B_2^N.
\end{equation}
Moreover,
\begin{equation}
\label{eq:poly-isometry-envelope}
  \norm{p_a}_{L^2(\rho)}=\norm{a}_2,
  \qquad \norm{E}_{L^2(\rho)}=\sqrt N,
  \qquad |p_a(x)|\leq\norm{a}_*E(x).
\end{equation}
\end{lemma}

\begin{proof}
Orthonormality gives
\[
\begin{aligned}
  \norm{p_a}_{L^2(\rho)}^2
    &= \sum_{i=1}^N a_i^2 = \norm{a}_2^2,\\
  \norm{E}_{L^2(\rho)}^2
    &= \sum_{i=1}^N \int q_i^2\,d\rho = N.
\end{aligned}
\]
By \eqref{eq:poly-V-star}, $|p_a(x)|\leq\norm{a}_*E(x)$ for every
$x\in\R^d$. Taking $L^2(\rho)$ norms gives the first inequality below.
The second follows from Cauchy--Schwarz and $V\subset S^{N-1}$:
\[
  \frac{\norm{a}_2}{\sqrt N}\leq\norm{a}_*\leq\norm{a}_2.
\]
The lower bound gives $\norm{a}_*>0$ for $a\neq0$. Absolute homogeneity
and the triangle inequality follow from the supremum in
\eqref{eq:poly-V-star}. Thus $\norm{\cdot}_*$ is a norm.
The two bounds also give \eqref{eq:poly-ball-inclusions}.
\end{proof}

\subsection{The complex evaluation cone}

Set $s=d+1$ and view $\R^N$ as the real subspace of $\C^N$.
Write $z=(z_0,z')\in\C^{d+1}$. For each basis polynomial, define
its degree-$m$ homogenization by
\[
  q_i(x)=\sum_{|\alpha|\leq m}c_{i,\alpha}x^\alpha,
  \qquad
  q_i^{\mathrm h}(z)
    =\sum_{|\alpha|\leq m}c_{i,\alpha}
      z_0^{m-|\alpha|}(z')^\alpha.
\]
Let $\nu_m(z)$ be the vector of all monomials of total degree $m$
in the coordinates of $z$, listed in a fixed order. Then
\[
  Q^{\mathrm h}(z)
  =\bigl(q_1^{\mathrm h}(z),\ldots,q_N^{\mathrm h}(z)\bigr)
  =T\nu_m(z)
\]
for a real $N\times N$ matrix $T$. Setting $z_0=1$ maps the monomials
in $\nu_m$ bijectively onto the monomial basis of $\mathcal P_m(\R^d)$.
So the $q_i^{\mathrm h}$ form a basis of the homogeneous polynomials
of degree $m$, and $T$ is invertible.

The map
\[
  F:\mathbf P^d(\C)\longrightarrow\mathbf P^{N-1}(\C),
  \qquad F([z])=[T\nu_m(z)]
\]
is the degree-$m$ Veronese embedding followed by the projective linear
automorphism induced by $T$. Its image $X$ is a compact complex submanifold with $\dim_{\C}X=d$. Define its affine cone by
\[
  \mathcal C=\{0\}\cup
  \{w\in\C^N\setminus\{0\}:[w]\in X\}.
\]
The map $w\mapsto[w]$ is continuous away from $0$. Since $X$ is closed
and $0\in\mathcal C$, the cone $\mathcal C$ is closed.

On the open set where $w_j\neq0$, the map $w\mapsto([w],w_j)$
has inverse $([u],t)\mapsto tu/u_j$. In these coordinates, the cone corresponds to the product of the open subset $\{[u]\in X:u_j\neq0\}$ and $\C\setminus\{0\}$. These coordinate domains cover $\mathcal C\setminus\{0\}$. Thus $\mathcal C\setminus\{0\}$ is a complex submanifold of complex
dimension $d+1=s$.

For $x\in\R^d$, we have $Q^{\mathrm h}(1,x)=Q(x)$, so
$[Q(x)]=F([1:x])\in X$. Hence $Q(x)/E(x)$ lies in
$\mathcal C\cap S^{2N-1}$. This intersection is closed, so
taking the closure gives
\begin{equation}
\label{eq:poly-V-cone}
  V\subseteq\mathcal C\cap S^{2N-1}.
\end{equation}

For a smooth real $k$-dimensional submanifold, let
$\operatorname{Vol}_k$ denote its induced Riemannian volume.
For a Borel set $A\subseteq\mathcal C$, we take $\operatorname{Vol}_{2s}(A)$ to be the Riemannian volume of
$A\setminus\{0\}$ in $\mathcal C\setminus\{0\}$. Let $\omega_k=\pi^{k/2}/\Gamma(k/2+1)$ be the volume of the Euclidean unit ball in $\R^k$.

\begin{lemma}[Volume of the complex cone]
\label{lem:poly-cone-volume}
For every $R>0$,
\begin{equation}
\label{eq:poly-cone-volume}
  \operatorname{Vol}_{2s}\bigl(\mathcal C\cap B(0,R)\bigr)
  =m^d\omega_{2s}R^{2s}.
\end{equation}
\end{lemma}

\begin{proof}
The Hopf map sends a unit vector $w$ to the complex line $[w]$.
The fiber through $w$ is the circle $\{e^{it}w:0\leq t<2\pi\}$.
Normalize the Fubini--Study metric so that the differential of the
Hopf map preserves the lengths of tangent vectors to the sphere
that are perpendicular to the fiber.
On the projective chart $z_j\neq0$, choose the representative $z$
with $z_j=1$. The K\"ahler form is
\[
  \omega_{\mathrm{FS}}
  =\frac{i}{2}\partial\bar\partial\log\norm{z}_2^2,
\]
where $i$ is the imaginary unit. 
In this normalization,
\[
  \int_{\mathbf P^d(\C)}\omega_{\mathrm{FS}}^d=\pi^d,
  \qquad
  \operatorname{Vol}_{2d}\bigl(\mathbf P^d(\C)\bigr)
  =\frac{\pi^d}{d!}.
\]
See \citet[Chapter~VI, (4.2) and Example~4.4]{Demailly2012},
whose Fubini--Study form is $\omega_{\mathrm{FS}}/\pi$.

Define
\[
  h([z])=\frac{\norm{T\nu_m(z)}_2^2}{\norm{z}_2^{2m}}.
\]
Under $z\mapsto\lambda z$ with $\lambda\neq0$, the numerator and
denominator both scale by $|\lambda|^{2m}$. The numerator is positive
for $z\neq0$ because $T$ is invertible and $\nu_m(z)\neq0$.
So $h$ is well defined and positive on $\mathbf P^d(\C)$. It is smooth in each projective chart, as the formula with $z_j=1$ shows. By the definition of $h$,
\[
  \log\norm{T\nu_m(z)}_2^2
  =m\log\norm{z}_2^2+\log h([z]).
\] 
The formula
for $\omega_{\mathrm{FS}}$ gives
\begin{equation}
\label{eq:poly-FS-pullback}
  F^*\omega_{\mathrm{FS}}
  =m\omega_{\mathrm{FS}}
   +\frac{i}{2}\partial\bar\partial\log h.
\end{equation}
The last term is exact:
\[
  \frac{i}{2}\partial\bar\partial\log h
  =\mathrm{d}\!\left[
    \frac{i}{4}(\bar\partial-\partial)\log h
  \right].
\]
In the expansion of $(F^*\omega_{\mathrm{FS}})^d$, every term
containing the exact form is itself exact, because the other factors
are closed. By Stokes' theorem, each such term integrates to zero
over $\mathbf P^d(\C)$. The K\"ahler volume formula gives
\begin{align}
\label{eq:poly-projective-volume}
  \operatorname{Vol}_{2d}(X)
  &=\frac{1}{d!}\int_{\mathbf P^d(\C)}
       (F^*\omega_{\mathrm{FS}})^d\notag\\
  &=\frac{m^d}{d!}\int_{\mathbf P^d(\C)}
       \omega_{\mathrm{FS}}^d
   =\frac{m^d\pi^d}{d!}.
\end{align}

Let $L=\mathcal C\cap S^{2N-1}$. The restricted Hopf map $L\to X$
is a Riemannian submersion with circle fibers of length $2\pi$.
Integrating over the fibers gives
\[
  \operatorname{Vol}_{2s-1}(L)
  =2\pi\operatorname{Vol}_{2d}(X).
\]
At radius $r$, the map $w\mapsto rw$ scales each of the $2s-1$
tangent directions of $L$ by $r$. The radial direction has unit
length and is orthogonal to these directions. Hence
\begin{align*}
  \operatorname{Vol}_{2s}\bigl(\mathcal C\cap B(0,R)\bigr)
  &=\operatorname{Vol}_{2s-1}(L)\int_0^R r^{2s-1}\,dr\\
  &=2\pi\,\frac{m^d\pi^d}{d!}\,\frac{R^{2s}}{2s}\\
  &=m^d\frac{\pi^s}{s!}R^{2s}
   =m^d\omega_{2s}R^{2s}.
\end{align*}
\end{proof}

\begin{lemma}[Local cone volume and evaluation nets]
\label{lem:poly-packing}
For every $v\in\mathcal C$ with $\norm{v}_2=1$ and every
$0<r\leq1/2$,
\begin{equation}
\label{eq:poly-local-volume}
  \operatorname{Vol}_{2s}\bigl(\mathcal C\cap B(v,r)\bigr)
  \geq\omega_{2s}r^{2s}.
\end{equation}
For every $0<\eta\leq1$, there is a Euclidean $\eta$-net for $V$
with centers in $V$ and cardinality at most
\begin{equation}
\label{eq:poly-net-cardinality}
  m^d\left(\frac{3}{\eta}\right)^{2s}.
\end{equation}
\end{lemma}

\begin{proof}
Set $M=\mathcal C\setminus\{0\}$ and $k=2s$.
Since $M$ is a complex submanifold of $\C^N$, its compact smooth
subdomains minimize volume among submanifolds with the same oriented
boundary \citep[Chapter~III, Theorem~1.26]{Demailly2012}.
The first variation formula then gives zero mean curvature on $M$.

Fix $v\in\mathcal C$ with $\norm{v}_2=1$. For $0<r<1$, the closed
ball $\overline{B(v,r)}$ avoids the origin. Write
\[
  R_v(w)=\norm{w-v}_2,
  \qquad
  A(r)=\operatorname{Vol}_k\bigl(M\cap B(v,r)\bigr).
\]

Let $(w-v)^\top$ be the orthogonal projection of $w-v$ onto $T_wM$.
Since $M$ has zero mean curvature, $\operatorname{div}_M((w-v)^\top)=k$.
Call $r$ a regular radius if $\nabla_M R_v(w)\neq0$ for every
$w\in M\cap\partial B(v,r)$. For such $r$, on $M\cap\partial B(v,r)$ we have
\[
  (w-v)^\top=r\nabla_M R_v,
  \qquad
  \nu=\frac{\nabla_M R_v}{\norm{\nabla_M R_v}_2},
\]
where $\nu$ is the outward unit normal to the boundary within $M$. The divergence theorem therefore gives
\begin{equation}
\label{eq:poly-divergence}
  kA(r)=r\int_{M\cap\partial B(v,r)}
               \norm{\nabla_M R_v}_2\,
               d\operatorname{Vol}_{k-1}.
\end{equation}

We next show that the critical set of $R_v^2$ on $M$ has zero
$k$-dimensional volume. For $u\in\mathcal C\cap S^{2N-1}$ and $t>0$,
\[
  \frac{\mathrm{d}}{\mathrm{d}t}R_v(tu)^2
  =2t+\norm{u-v}_2^2-2.
\]
Since the radial direction is tangent to $M$, the radial derivative
vanishes at every critical point of $R_v^2$ on $M$. This derivative
vanishes at most once on each ray, so integration in radial coordinates
shows that the critical set has zero $k$-dimensional volume. On
$M\setminus\{v\}$, the identity $\nabla_M(R_v^2)=2R_v\nabla_M R_v$
shows that $R_v$ and $R_v^2$ have the same critical points. Since the critical set has zero volume, the coarea formula shows
that $A$ is locally absolutely continuous on $(0,1)$ and that,
for almost every $r\in(0,1)$,
\begin{equation}
\label{eq:poly-coarea}
  A'(r)=\int_{M\cap\partial B(v,r)}
          \frac{1}{\norm{\nabla_M R_v}_2}\,
          d\operatorname{Vol}_{k-1}.
\end{equation}
By Sard's theorem applied to $R_v$ on $M\setminus\{v\}$, almost
every $r\in(0,1)$ is regular. Since $\norm{\nabla_M R_v}_2\leq1$,
\eqref{eq:poly-divergence} and \eqref{eq:poly-coarea} give
\[
  \left(\frac{A(r)}{r^k}\right)'
  =\frac{rA'(r)-kA(r)}{r^{k+1}}\geq0
  \qquad\text{for almost every }r\in(0,1).
\]
Because $A$ is locally absolutely continuous on $(0,1)$, so is
$A(r)/r^k$. The derivative inequality therefore shows that
$A(r)/r^k$ is nondecreasing. Since $M$ is smooth at $v$,
$\lim_{r\downarrow0}A(r)/r^k=\omega_k$.
This proves \eqref{eq:poly-local-volume}.

By compactness, $V$ has a finite maximal $\eta$-separated subset
$\{v_1,\ldots,v_J\}$. It is an $\eta$-net for $V$ by maximality.
The open balls $B(v_j,\eta/2)$ are pairwise disjoint and lie in
$B(0,1+\eta/2)$.
By \eqref{eq:poly-V-cone}, \eqref{eq:poly-local-volume}, and
\eqref{eq:poly-cone-volume},
\[
\begin{aligned}
  J\omega_{2s}(\eta/2)^{2s}
  &\leq \sum_{j=1}^J
    \operatorname{Vol}_{2s}\bigl(\mathcal C\cap B(v_j,\eta/2)\bigr)\\
  &\leq \operatorname{Vol}_{2s}
    \bigl(\mathcal C\cap B(0,1+\eta/2)\bigr)\\
  &=m^d\omega_{2s}(1+\eta/2)^{2s}.
\end{aligned}
\]
Since $0<\eta\leq1$, it follows that
\[
  J\leq m^d\left(\frac{2+\eta}{\eta}\right)^{2s}
   \leq m^d\left(\frac{3}{\eta}\right)^{2s}.
\]
This proves \eqref{eq:poly-net-cardinality}.
\end{proof}

\subsection{Gaussian small balls and coefficient volume}

We use the $2^{-j}$-nets from Lemma~\ref{lem:poly-packing}
in a chaining argument following
\citet[Section~4, equation~(4.12)]{Talagrand1993}
(see also \citealp[Section~2, Lemma~13]{AurzadaLifshits2008}).
The factor $m^d$ in the bound on the net sizes yields a
small-ball probability lower bound of the form $\exp(-c_dm^d)$.

\begin{lemma}[A Gaussian small-ball estimate]
\label{lem:poly-small-ball}
There are positive constants $T_d$ and $c_d$, depending only on $d$,
such that a standard Gaussian vector $g$ in $\R^N$ satisfies
\begin{equation}
\label{eq:poly-small-ball}
  \Pp\left\{\sup_{v\in V}|\langle g,v\rangle|\leq T_d\right\}
  \geq \exp(-c_dm^d).
\end{equation}
\end{lemma}

\begin{proof}
For each integer $j\geq0$, choose a $2^{-j}$-net $A_j\subset V$
using Lemma~\ref{lem:poly-packing}. Then
\begin{equation}
\label{eq:poly-dyadic-cardinality}
  |A_j|\leq m^d3^{2s}2^{2sj}.
\end{equation}
For $j\geq1$ and $a\in A_j$, choose $\pi_j(a)\in A_{j-1}$ with
$\norm{a-\pi_j(a)}_2\leq2^{1-j}$. Define
\[
  u_0(a)=a\quad(a\in A_0),
  \qquad
  u_j(a)=a-\pi_j(a)\quad(j\geq1,\ a\in A_j).
\]
For $a\in A_0$, we have $\norm{u_0(a)}_2=1$. Together with the
bound for $j\geq1$, this gives $\norm{u_j(a)}_2\leq2^{1-j}$
for all $j\geq0$ and $a\in A_j$.

Set $b=2s\log 2+2$ and $K=\sqrt{8b}$. For an integer $J\geq0$,
let $\mathcal E_J$ be the event that
\[
  |\langle g,u_j(a)\rangle|
  \leq K\sqrt{j+1}\,2^{-j}
  \qquad\text{for all }0\leq j\leq J\text{ and }a\in A_j.
\]
Since $\operatorname{Var}(\langle g,u_j(a)\rangle)
=\norm{u_j(a)}_2^2\leq4\cdot4^{-j}$, the Gaussian tail bound gives
\[
  \Pp\left\{|\langle g,u_j(a)\rangle|>
      K\sqrt{j+1}\,2^{-j}\right\}
  \leq 2e^{-b(j+1)}=:q_j\leq\tfrac12.
\]
Each zero-variance variable is zero almost surely, so its bound
holds automatically. For the remaining variables, \v{S}id\'ak's
inequality \citep[Corollary~1]{Sidak1967} applies even when their
covariance matrix is singular. Using the marginal probability
bounds $1-q_j$, we obtain
\[
  \Pp(\mathcal E_J)\geq\prod_{j=0}^{J}(1-q_j)^{|A_j|}.
\]
Using $\log(1-q)\geq-2q$ for $0\leq q\leq1/2$ and
\eqref{eq:poly-dyadic-cardinality}, we get
\begin{align*}
  \log\Pp(\mathcal E_J)
  &\geq -2\sum_{j=0}^{J}|A_j|q_j\\
  &\geq -4m^d3^{2s}e^{-b}\sum_{j=0}^{\infty}e^{-2j}
   =-c_dm^d,
  \qquad c_d=\frac{4\cdot3^{2s}e^{-b}}{1-e^{-2}}.
\end{align*}
Since $\mathcal E_{J+1}\subseteq\mathcal E_J$ for every $J\geq0$,
\[
  \Pp\left(\bigcap_{J=0}^{\infty}\mathcal E_J\right)
  =\lim_{J\to\infty}\Pp(\mathcal E_J)
  \geq e^{-c_dm^d}.
\]

On this intersection, fix any integer $J\geq0$ and any $a\in A_J$.
Set $a_J=a$ and $a_{j-1}=\pi_j(a_j)$ for $j=J,\ldots,1$.
Telescoping gives
\[
  |\langle g,a\rangle|
  =\left|\sum_{j=0}^{J}\langle g,u_j(a_j)\rangle\right|
  \leq K\sum_{j=0}^{J}\sqrt{j+1}\,2^{-j}
  \leq T_d,
  \qquad
  T_d=K\sum_{j=0}^{\infty}\sqrt{j+1}\,2^{-j}<\infty.
\]
For every $v\in V$, choose $a_J\in A_J$ with
$\norm{v-a_J}_2\leq2^{-J}$. Then
\[
  |\langle g,v\rangle|
  \leq |\langle g,a_J\rangle|+\norm{g}_2\norm{v-a_J}_2
  \leq T_d+2^{-J}\norm{g}_2.
\]
Letting $J\to\infty$ gives $\sup_{v\in V}|\langle g,v\rangle|\leq T_d$
on this intersection. The probability bound above then gives
\eqref{eq:poly-small-ball}.
\end{proof}

\begin{lemma}[Volume of the evaluation-norm unit ball]
\label{lem:poly-body-volume}
There is a constant $a_d>0$, depending only on $d$, such that
\begin{equation}
\label{eq:poly-body-volume}
  \operatorname{Vol}_{N}(B_*)\geq a_d^N.
\end{equation}
\end{lemma}

\begin{proof}
The event in \eqref{eq:poly-small-ball} is
$\{\norm{g}_*\leq T_d\}=\{g\in T_dB_*\}$.
The standard Gaussian density in $\R^N$ is at most $(2\pi)^{-N/2}$,
so
\begin{align*}
  e^{-c_dm^d}
  &\leq \Pp\{g\in T_dB_*\}\\
  &\leq (2\pi)^{-N/2}\operatorname{Vol}_{N}(T_dB_*)
   = (2\pi)^{-N/2}T_d^N\operatorname{Vol}_{N}(B_*).
\end{align*}
Since $m^d\leq(m+1)\cdots(m+d)=d!N$, it follows that
\[
  \operatorname{Vol}_{N}(B_*)
  \geq \left(\frac{\sqrt{2\pi}}{T_d}\right)^N e^{-c_dm^d}
  \geq \left(\frac{\sqrt{2\pi}}{T_d}e^{-c_d d!}\right)^N.
\]
This proves \eqref{eq:poly-body-volume} with
$a_d=(\sqrt{2\pi}/T_d)e^{-c_d\,d!}$, which depends only on $d$.
\end{proof}

\subsection{Construction of the brackets}

\begin{proof}[Proof of Theorem~\ref{thm:poly}]
Fix $0<\varepsilon\leq1$. First suppose $m\geq1$ and set
$\delta=\varepsilon/(2\sqrt N)$. Choose a maximal $\delta$-separated subset $\{a_1,\ldots,a_J\}$
of $B_2^N$ with respect to $\norm{\cdot}_*$. By maximality, it is a $\delta$-net for $B_2^N$ in this norm.

The translated balls $a_j+(\delta/2)B_*$ have disjoint interiors.
By \eqref{eq:poly-ball-inclusions}, each lies in
\[
  B_2^N+(\delta/2)B_*
  \subset \left(1+\frac{\delta\sqrt N}{2}\right)B_2^N
  =\left(1+\frac{\varepsilon}{4}\right)B_2^N.
\]
Comparing volumes gives
\begin{equation}
\label{eq:poly-volumetric-cover}
  J\leq
  \left(\frac{4\sqrt N}{\varepsilon}\right)^N
  \left(1+\frac{\varepsilon}{4}\right)^N
  \frac{\operatorname{Vol}_{N}(B_2^N)}
       {\operatorname{Vol}_{N}(B_*)}.
\end{equation}
By Stirling's inequality, for an absolute constant $C$,
\[
  \operatorname{Vol}_{N}(B_2^N)
  =\frac{\pi^{N/2}}{\Gamma(N/2+1)}
  \leq\left(\frac{C}{\sqrt N}\right)^N.
\]
Lemma~\ref{lem:poly-body-volume} gives
$\operatorname{Vol}_{N}(B_*)\geq a_d^N$. Substituting the two volume
bounds into \eqref{eq:poly-volumetric-cover} gives
\begin{equation}
\label{eq:poly-cover-count}
  J\leq
  \left(\frac{4C(1+\varepsilon/4)}{a_d\varepsilon}\right)^N
  \leq\left(\frac{C_d}{\varepsilon}\right)^N.
\end{equation}

For each $p_a\in\mathcal F_m(\rho)$, the isometry in
\eqref{eq:poly-isometry-envelope} gives $a\in B_2^N$.
Choose $a_j$ with $\norm{a-a_j}_*\leq\delta$. Then
\[
  |p_a(x)-p_{a_j}(x)|
  =|p_{a-a_j}(x)|
  \leq\norm{a-a_j}_* E(x)
  \leq\delta E(x),
  \qquad x\in\R^d.
\]
So the $J$ brackets
\[
  [p_{a_j}-\delta E,\ p_{a_j}+\delta E],
  \qquad 1\leq j\leq J,
\]
cover $\mathcal F_m(\rho)$. Each has $L^2(\rho)$ width
$2\delta\norm{E}_{L^2(\rho)}=2\delta\sqrt N=\varepsilon$.
Taking logarithms in \eqref{eq:poly-cover-count} gives
\eqref{eq:poly-entropy} for $m\geq1$.

For $m=0$, we have $N=1$ and
$\mathcal F_0(\rho)=\{a q_1:-1\leq a\leq1\}$, where
$q_1=\rho(\R^d)^{-1/2}$. Divide $[-1,1]$ into
$\lceil2/\varepsilon\rceil$ intervals of length at most $\varepsilon$.
Each interval $[a,b]$ gives a bracket $[a q_1,b q_1]$ of width
$b-a$. Since $\lceil2/\varepsilon\rceil\leq3/\varepsilon$, the upper
bound holds also for $m=0$ after enlarging $C_d$ if needed.

For the lower bound, consider any cover of $\mathcal F_m(\rho)$
by $J$ brackets of width at most $\varepsilon$.
From each bracket $[\ell,u]$ meeting $\mathcal F_m(\rho)$,
choose $p_a\in\mathcal F_m(\rho)$. For any
$p_b\in\mathcal F_m(\rho)$ in the same bracket,
we have $|p_b-p_a|\leq u-\ell$ pointwise. Hence
\[
  \norm{b-a}_2
  =\norm{p_b-p_a}_{L^2(\rho)}
  \leq\norm{u-\ell}_{L^2(\rho)}
  \leq\varepsilon.
\]
Thus at most $J$ Euclidean balls of radius $\varepsilon$ cover
$B_2^N$. Comparing volumes gives
\[
  \operatorname{Vol}_{N}(B_2^N)
  \leq J\varepsilon^N\operatorname{Vol}_{N}(B_2^N),
  \qquad J\geq\varepsilon^{-N}.
\]
Taking logarithms proves the lower bound.
\end{proof}

\section{Local bracketing of Gaussian mixture densities}
\label{sec:local}

Write $h_G=\sqrt{f_G}$. For any two mixture densities,
\[
  H^2(f_G,f_{G_0})
  =\norm{h_G-h_{G_0}}_{L^2(\R^d)}^2\leq2.
\]
Let $K_r=[-r,r]^d$ for $r>0$. For $R\geq1$, define
\[
  \mathcal F_R=\{f_G:G\in\mathcal P(K_R)\},
  \qquad
  \mathcal H_R=\{h_G:G\in\mathcal P(K_R)\}.
\]
Constants indexed by $d$ depend only on the dimension.

We adapt the weighted decomposition in \citet[Section~3.2]{HM2026}
and divide $h_G-h_{G_0}$ by the local radius $s$. On a larger cube
in the observation space, we use a Taylor polynomial to approximate
the factor multiplying the weight. On the same cube, we use a grid
to cover the weights at a fixed multiplicative accuracy. For an
analogous cover for univariate finite Gaussian mixtures, see
\citet[Proposition~3.8]{HM2026}. Combining these approximations with Theorem~\ref{thm:poly}
and a tail bound outside the cube, we construct $L^2(\R^d)$
brackets of width at most $s/2$ for the radius-$s$ ball in
$\mathcal H_R$ centered at $h_{G_0}$.

\begin{proposition}[Bracketing at half the radius]
\label{prop:local-half}
There is a constant $C_d<\infty$ such that for every $R\geq1$,
$0<s\leq\sqrt2$, and $h_0\in\mathcal H_R$,
\begin{equation}
\label{eq:local-half}
  \log N_{[\,]}\!\left(
    \frac{s}{2},
    \{h\in\mathcal H_R:\norm{h-h_0}_2\leq s\},
    L^2(\R^d)
  \right)
  \leq C_d\bigl[1+R^2+\log_+(1/s)\bigr]^d.
\end{equation}
The bracket endpoints may be taken nonnegative.
\end{proposition}

\begin{proof}
Choose $G_0\in\mathcal P(K_R)$ with $h_0=h_{G_0}$. Set
\[
  L=1+R^2+\log_+(1/s),\qquad
  T=A_d\sqrt L,\qquad I=K_T,\qquad m=\lceil B_dL\rceil.
\]
We choose $A_d\geq2$ first and then $B_d>0$. Both constants depend
only on $d$.

For $G\in\mathcal P(K_R)$, define
\begin{align*}
  \psi_G(x)
  &=\frac{f_G(x)}{\phi_d(x)}
   =\int_{K_R} e^{x\cdot\theta-\norm{\theta}_2^2/2}\,dG(\theta),\\
  W_G(x)
  &=\bigl(\sqrt{\psi_G(x)}+\sqrt{\psi_{G_0}(x)}\bigr)^2,
  \qquad
  V_G(x)=\frac{\sqrt{\phi_d(x)}}{\sqrt{W_G(x)}}.
\end{align*}
For each $G$ with $\norm{h_G-h_0}_2\leq s$, set
\begin{equation}
\label{eq:normalized-local}
  z_G=\frac{h_G-h_0}{s}
      =\frac{\psi_G-\psi_{G_0}}{s}\,V_G,
  \qquad \norm{z_G}_2\leq1.
\end{equation}
We use the same radius $s$ in the denominator for every $h_G$ in the ball.

Fix $\eta=1/128$. For $x\in I^c$, both $h_G(x)$ and $h_{G_0}(x)$
are at most $\sup_{\theta\in K_R}\sqrt{\phi_d(x-\theta)}$. So
\[
  |z_G(x)|\leq F(x),
  \qquad
  F(x)=\frac{2}{s}\sup_{\theta\in K_R}\sqrt{\phi_d(x-\theta)}.
\]
Minimizing $\norm{x-\theta}_2^2$ over the cube gives
\[
  \sup_{\theta\in K_R}\phi_d(x-\theta)
  =(2\pi)^{-d/2}\prod_{j=1}^d
    \exp\!\left\{-\frac{((|x_j|-R)_+)^2}{2}\right\}.
\]
The one-dimensional factor integrates to $2R+\sqrt{2\pi}$ over
$\R$. Since $T-R\geq(A_d-1)\sqrt L\geq1$, its integral over
$|t|>T$ satisfies
\[
  \int_{|t|>T}\exp\!\left\{-\frac{((|t|-R)_+)^2}{2}\right\}\,dt
  =2\int_{T-R}^\infty e^{-u^2/2}\,du
  \leq2e^{-(T-R)^2/2}.
\]
Since $I^c$ is the union of the $d$ sets $\{|x_j|>T\}$, we obtain
\begin{equation}
\label{eq:tail-envelope}
  \norm{F1_{I^c}}_2^2
  \leq C_d s^{-2}(1+R)^{d-1}e^{-(T-R)^2/2}.
\end{equation}
Here $1_A$ denotes the indicator of $A$. Using $s^{-2}\leq e^{2L}$,
$\log(1+R)\leq L$, and $T-R\geq(A_d-1)\sqrt L$, we bound the
right-hand side of \eqref{eq:tail-envelope} by
\[
  C_d\exp\!\left\{
    \left[d+1-\frac{(A_d-1)^2}{2}\right]L
  \right\}.
\]
Since $L\geq2$, we can choose $A_d$ large enough, depending only on
$d$, to ensure $\norm{F1_{I^c}}_2\leq\eta$.

On $I$, approximate $\psi_G$ by the polynomial of total degree at
most $m$ given by
\[
  p_G(x)=\int_{K_R} e^{-\norm{\theta}_2^2/2}
            \sum_{j=0}^m\frac{(x\cdot\theta)^j}{j!}\,dG(\theta).
\]
Set $a=dRT$. Since $|x\cdot\theta|\leq a$ for $x\in I$ and
$\theta\in K_R$, Taylor's formula gives
\[
  |\psi_G(x)-p_G(x)|\leq e^a\frac{a^{m+1}}{(m+1)!},
  \qquad x\in I.
\]
For every $\theta\in K_R$, we have $\norm{\theta}_2^2\leq dR^2$.
Hence $\psi_G$ and $\psi_{G_0}$ are both at least
$e^{-a-dR^2/2}$ on $I$, so
\[
  V_G(x)\leq\frac12\sqrt{\phi_d(x)}\,e^{a/2+dR^2/4},
  \qquad x\in I.
\]
With $q_G=(p_G-p_{G_0})/s$, equation~\eqref{eq:normalized-local}
gives, for $x\in I$,
\begin{align}
  |z_G(x)-q_G(x)V_G(x)|
  &\leq s^{-1}V_G(x)
     \bigl(|\psi_G(x)-p_G(x)|+|\psi_{G_0}(x)-p_{G_0}(x)|\bigr)
     \notag\\
  &\leq s^{-1}\sqrt{\phi_d(x)}\,e^{3a/2+dR^2/4}
          \frac{a^{m+1}}{(m+1)!}
   =:E(x).
  \label{eq:poly-error-envelope}
\end{align}
Choose $B_d>edA_d$ large enough that
\[
  2B_d\log\!\left(\frac{B_d}{edA_d}\right)
  \geq 2+3dA_d+d/2+\log(1/\eta).
\]
Using $a\leq dA_dL$, $R^2\leq L$, $s^{-2}\leq e^{2L}$, and
$(m+1)!\geq((m+1)/e)^{m+1}$, we get
\begin{align*}
  \norm{E}_{L^2(I)}^2
  &\leq e^{(2+3dA_d+d/2)L}
       \left(\frac{edA_dL}{m+1}\right)^{2(m+1)}\\
  &\leq e^{(2+3dA_d+d/2)L}
       \left(\frac{edA_d}{B_d}\right)^{2B_dL}\\
  &\leq \eta^L\leq\eta^2.
\end{align*}
The second inequality uses $m+1\geq B_dL$ and $edA_d/B_d<1$. Therefore
\begin{equation}
\label{eq:weighted-poly-norm}
  \norm{q_GV_G}_{L^2(I)}
  \leq\norm{z_G}_{L^2(I)}+\norm{E}_{L^2(I)}
  \leq1+\eta.
\end{equation}

Let $g_G=\log W_G$, so that $V_G=\sqrt{\phi_d}\,e^{-g_G/2}$.
Differentiating under the integral gives
\[
  \partial_i\log\psi_G(x)
  =\frac{\int_{K_R}\theta_i
       e^{x\cdot\theta-\norm{\theta}_2^2/2}\,dG(\theta)}
       {\psi_G(x)}
  \in[-R,R].
\]
The chain rule gives
\[
  \partial_i g_G
  =\frac{\sqrt{\psi_G}\,\partial_i\log\psi_G
         +\sqrt{\psi_{G_0}}\,\partial_i\log\psi_{G_0}}
        {\sqrt{\psi_G}+\sqrt{\psi_{G_0}}}.
\]
This is a weighted average of two numbers in $[-R,R]$.
At $x=0$, both $\psi_G(0)$ and $\psi_{G_0}(0)$ lie in
$[e^{-dR^2/2},1]$. Thus
\begin{equation}
\label{eq:weight-lipschitz}
  |\partial_i g_G|\leq R,\qquad
  \log 4-\frac{dR^2}{2}\leq g_G(0)\leq\log 4.
\end{equation}

Fix $\kappa>0$. Set $M=\lceil4dRT/\kappa\rceil$ and use the grid
$(T/M)\{-M,\ldots,M\}^d$ on $I$. At each grid point $x$, record
the integer $\lfloor4g_G(x)/\kappa\rfloor$. By
\eqref{eq:weight-lipschitz}, there are at most
$C_{d,\kappa}(1+R^2)$ possible integers at $0$. For adjacent grid
points $x,x'$, we have
\[
  |g_G(x)-g_G(x')|\leq RT/M\leq\kappa/(4d).
\]
Their recorded integers therefore differ by at most one. Order the grid points $(T/M)k$ by increasing $\sum_{i=1}^d|k_i|$,
breaking ties arbitrarily. For each point other than $0$, fix one
adjacent point earlier in the order. Given the record at that
neighbor, there are at most three choices for the new record. The logarithm
of the number of arrays is therefore at most
\[
  \log\!\bigl(C_{d,\kappa}(1+R^2)\bigr)
  +\bigl((2M+1)^d-1\bigr)\log 3
  \leq C_{d,\kappa}L^d,
\]
where the last inequality uses $RT\leq A_dL$.

If two log-weight functions $g$ and $\widetilde g$ give the same
array, their values differ by at most $\kappa/4$ at each grid point. For $x\in I$, choose a grid
point $x'$ within $T/M$ in each coordinate. Then
\[
  |g(x)-\widetilde g(x)|
  \leq 2dRT/M+\kappa/4
  \leq 3\kappa/4
  \leq\kappa.
\]
For each array that occurs, choose one function $g_j$ that gives
that array and set $V_j=\sqrt{\phi_d}\,e^{-g_j/2}$ on $I$. If $g_G$ gives the same
array as $g_j$, then on $I$,
\begin{equation}
\label{eq:weight-net-ratio}
  e^{-\kappa/2}
  \leq\frac{V_G}{V_j}
  =e^{-(g_G-g_j)/2}
  \leq e^{\kappa/2}.
\end{equation}

The measure $\rho_j(dx)=V_j(x)^2\,1_I(x)\,dx$ has finite moments
of every order because $I$ is compact and $V_j$ is continuous on
$I$. Its density is positive on the interior of $I$, so every
nonzero polynomial has positive $L^2(\rho_j)$ norm. By
\eqref{eq:weight-net-ratio} and \eqref{eq:weighted-poly-norm},
\[
  \norm{q_G}_{L^2(\rho_j)}
  =\norm{q_GV_j}_{L^2(I)}
  \leq e^{\kappa/2}\norm{q_GV_G}_{L^2(I)}
  \leq e^{\kappa/2}(1+\eta)\leq2
\]
for sufficiently small $\kappa$. Fix $0<\eta_p\leq1$. Apply Theorem~\ref{thm:poly} at width
$\eta_p/2$ to the $L^2(\rho_j)$ unit ball of polynomials of total
degree at most $m$. Multiplying the bracket endpoints by $2$
gives brackets of $L^2(\rho_j)$ width at most $\eta_p$ that cover
all $q_G$ corresponding to $j$. The logarithm of the number of
brackets is at most
\[
  C_{d,\eta_p}\binom{m+d}{d}\leq C_{d,\eta_p}L^d.
\]

For each $j$, retain the brackets $[\ell,u]$ that contain a polynomial
$q$ with $\norm{q}_{L^2(\rho_j)}\leq2$. For such a bracket,
\[
  \max\{\norm{\ell}_{L^2(\rho_j)},\norm{u}_{L^2(\rho_j)}\}
  \leq \norm{q}_{L^2(\rho_j)}+\norm{u-\ell}_{L^2(\rho_j)}
  \leq 2+\eta_p.
\]
Set $\beta=e^{\kappa/2}-1$ and $B=|\ell|+|u|$. If $G$ is assigned to $j$ and $\ell\leq q_G\leq u$, then
\eqref{eq:weight-net-ratio} gives $|V_G/V_j-1|\leq\beta$ on $I$.
Hence
\[
  \left|q_G\frac{V_G}{V_j}-q_G\right|
  \leq \beta|q_G|\leq\beta B.
\]
Thus
\[
  V_j(\ell-\beta B)\leq q_GV_G\leq V_j(u+\beta B).
\]
The $L^2(I)$ width of this bracket is bounded by
\begin{align*}
  \norm{V_j(u-\ell+2\beta B)}_{L^2(I)}
  &=\norm{u-\ell+2\beta B}_{L^2(\rho_j)}\\
  &\leq\eta_p+4\beta(2+\eta_p).
\end{align*}
Choose fixed $\eta_p,\kappa>0$ small enough that
$\eta_p+4\beta(2+\eta_p)\leq1/8$ and
$e^{\kappa/2}(1+\eta)\leq2$.

On $I$, subtract $E$ from the lower endpoint and add $E$ to the
upper endpoint. On $I^c$, use $[-F,F]$. By
\eqref{eq:tail-envelope} and \eqref{eq:poly-error-envelope}, the
resulting bracket $[\underline z,\overline z]$ contains $z_G$ and has
$L^2(\R^d)$ width at most
\[
  \frac18+2\norm{E}_{L^2(I)}+2\norm{F1_{I^c}}_2
  \leq \frac18+4\eta<\frac12.
\]
Since $h_G=h_0+sz_G\geq0$, the bracket
\[
  \bigl[(h_0+s\underline z)_+,\,(h_0+s\overline z)_+\bigr]
\]
contains $h_G$. Its width is at most $s/2$ because the positive-part
map is increasing and $1$-Lipschitz. There are at most
$\exp(C_{d,\kappa}L^d)$ weight representatives and at most
$\exp(C_{d,\eta_p}L^d)$ polynomial brackets for each representative.
Since $\kappa$ and $\eta_p$ are fixed absolute constants, the total
number is at most $\exp(C_dL^d)$. This proves \eqref{eq:local-half}.
\end{proof}

\begin{corollary}[Local entropy integral]
\label{cor:local-integral}
Uniformly over $R\geq1$, $h_0\in\mathcal H_R$, and
$0<u\leq\delta\leq\sqrt2$,
\begin{equation}
\label{eq:recursive-entropy}
  \log N_{[\,]}\!\left(
    u,\{h\in\mathcal H_R:\norm{h-h_0}_2\leq\delta\},L^2(\R^d)
  \right)
  \leq C_d[1+R^2+\log_+(1/u)]^d[1+\log(\delta/u)].
\end{equation}
Define
\[
  \mathcal J_{[\,]}(\mathcal F_R;\delta)
  =
  \sup_{f_0\in\mathcal F_R}
  \int_0^\delta
    \sqrt{1+\log N_{[\,]}\!\left(
      u,\{f\in\mathcal F_R:H(f,f_0)\leq\delta\},H
    \right)}
  \,du.
\]
Then, for $0<\delta\leq\sqrt2$,
\begin{equation}
\label{eq:entropy-integral}
  \mathcal J_{[\,]}(\mathcal F_R;\delta)
  \leq C_d\delta[1+R^2+\log_+(1/\delta)]^{d/2}.
\end{equation}
\end{corollary}

\begin{proof}
We iterate Proposition~\ref{prop:local-half} as in
\citet[Appendix~B.6, Step~2]{HM2026}. For a related local bracketing
argument based on a decomposition into shells, see \citet[Theorem~2.2]{GVH2014}.
Set $s_j=\delta 2^{-j}$. Apply the proposition first to the ball of
radius $s_0$ centered at $h_0$. At each subsequent stage $j\geq1$, for each bracket of width at most
$s_j$ that meets $\mathcal H_R$, choose a member of $\mathcal H_R$
in the bracket as the next center. Every member of $\mathcal H_R$ in
that bracket lies within $L^2$ distance $s_j$ of the center. Apply
the proposition to this ball to obtain brackets of width at most $s_{j+1}$.

Let $k=\max\{1,\lceil\log_2(\delta/u)\rceil\}$. Then $s_k\leq u$ and
$s_j\geq u$ for $0\leq j<k$. The logarithm of the number of brackets
after $k$ steps is at most
\[
  C_d\sum_{j=0}^{k-1}[1+R^2+\log_+(1/s_j)]^d
  \leq C_d[1+R^2+\log_+(1/u)]^d[1+\log(\delta/u)].
\]
This proves \eqref{eq:recursive-entropy}.

Let $[a,b]$ be one of the nonnegative brackets that contains $h_G$.
Since $0\leq a\leq h_G\leq b\leq h_G+(b-a)$, both endpoints lie
in $L^2(\R^d)$. Thus $[a^2,b^2]$ is a bracket with integrable
endpoints containing $f_G$. Its Hellinger width is
$\norm{\sqrt{b^2}-\sqrt{a^2}}_2=\norm{b-a}_2$.

In the entropy integral, set $u=\delta e^{-t}$ and
$A=1+R^2+\log_+(1/\delta)$. Then $\log(\delta/u)=t$ and
$\log_+(1/u)\leq\log_+(1/\delta)+t$. By
\eqref{eq:recursive-entropy} and $A+t\leq A(1+t)$,
\begin{align*}
  \mathcal J_{[\,]}(\mathcal F_R;\delta)
  &\leq C_d\delta\int_0^\infty
    e^{-t}(A+t)^{d/2}\sqrt{1+t}\,dt\\
  &\leq C_d\delta A^{d/2}\int_0^\infty
    e^{-t}(1+t)^{(d+1)/2}\,dt\\
  &\leq C_d\delta A^{d/2}.
\end{align*}
This proves \eqref{eq:entropy-integral}.
\end{proof}

\section{Likelihood risk and truncation}
\label{sec:risk}

\begin{theorem}[Risk on a compact mixing class]
\label{thm:compact}
For every $d\geq1$, there is a constant $C_d<\infty$ such that, for $R\geq1$,
$n\geq2$, and any measurable maximum likelihood estimator
$\widehat f_{n,R}$ over $\mathcal F_R$,
\begin{equation}
\label{eq:compact-risk}
 \sup_{G\in\mathcal P(K_R)}
 \E_{f_G^{\otimes n}}H^2(\widehat f_{n,R},f_G)
 \leq C_d\min\!\left\{1,\frac{(1+R^2+\log n)^d}{n}\right\}.
\end{equation}
Since $C_d$ depends only on $d$, the bound also holds with $R=R_n$ for any deterministic sequence $R_n\geq1$.
\end{theorem}

\begin{proof}
Suppose $G_j\Rightarrow G$ on $K_R$. Since the Gaussian kernel is
bounded and continuous,
\[
  f_{G_j}(x)\to f_G(x)>0,\qquad
  \log f_{G_j}(x)\to\log f_G(x)
  \quad\text{for every }x\in\R^d.
\]
By Scheff\'e's lemma,
\[
  H^2(f_{G_j},f_G)\leq\norm{f_{G_j}-f_G}_1\to0.
\]
Thus, for each $G_0\in\mathcal P(K_R)$ and $\delta>0$, the set
\[
  \{G\in\mathcal P(K_R):H(f_G,f_{G_0})\leq\delta\}
\]
is closed in the weak topology. This set and $\mathcal P(K_R)$
are both compact and metrizable. A countable weakly dense subset
of each gives a pointwise dense class of log likelihood functions.

Write $f_0=f_{G_0}$. For every $G\in\mathcal P(K_R)$,
\[
  \frac{f_G+f_0}{2}=f_{(G+G_0)/2}\in\mathcal F_R.
\]
So every midpoint density within Hellinger distance $\delta$ of
$f_0$ lies in the radius-$\delta$ ball in $\mathcal F_R$ centered at
$f_0$. The entropy bound \eqref{eq:entropy-integral} therefore
controls the local midpoint bracketing integral in
\citet[Theorem~7.2.1 and p.~550]{GN2016}. For an earlier likelihood-ratio
bound under a local bracketing entropy condition, see
\citet[Theorem~1 and the remark on p.~349]{WS1995}.

We apply the MLE risk bound of \citet[Lemma~2.1]{HM2026}.
Since $\sqrt{1+x}\geq\max\{1,\sqrt{x}\}$ for $x\geq0$,
our entropy integral dominates the one in that lemma.
By \eqref{eq:entropy-integral}, we may take
\[
  \Psi(\delta)=C_d\delta[1+R^2+\log_+(1/\delta)]^{d/2},
  \qquad 0<\delta\leq\sqrt2.
\]
The ratio
\[
  \frac{\Psi(\delta)}{\delta^2}
  =\frac{C_d}{\delta}[1+R^2+\log_+(1/\delta)]^{d/2}
\]
is nonincreasing. For a universal constant $C$, the likelihood bound gives expected squared Hellinger loss at most $2\delta_n^2$ if
$0<\delta_n\leq\sqrt2$ and $\sqrt n\,\delta_n^2\geq C\Psi(\delta_n)$.

Set $A_n=1+R^2+\log n$ and $\delta_n^2=K_d A_n^d/n$, where
$K_d\geq1$ and $\sqrt{K_d}\geq CC_d$. If $\delta_n\leq1/2$, then
$\delta_n\geq n^{-1/2}$ gives $\log(1/\delta_n)\leq(\log n)/2$.
Consequently,
\begin{align*}
  \sqrt n\,\delta_n^2
  &=\sqrt{K_d}\,\delta_n A_n^{d/2}\\
  &\geq CC_d\,\delta_n[1+R^2+\log(1/\delta_n)]^{d/2}\\
  &=C\Psi(\delta_n).
\end{align*}
The likelihood bound therefore gives risk at most $2\delta_n^2$.
If $\delta_n>1/2$, then $H^2\leq2<8\delta_n^2$. Combining both cases with $H^2\leq2$ gives
\[
  \sup_{G\in\mathcal P(K_R)}
  \E_{f_G^{\otimes n}}H^2(\widehat f_{n,R},f_G)
  \leq \min\{2,8\delta_n^2\}
  \leq 8K_d\min\!\left\{1,\frac{A_n^d}{n}\right\}.
\]
This proves \eqref{eq:compact-risk}.
\end{proof}

If all observations lie in $K_R$, then their convex hull lies in
$K_R$ by convexity. By Lemma~\ref{lem:convex-hull-support}, every
unrestricted maximizer then belongs to $\mathcal P(K_R)$ and hence
also maximizes the likelihood over this class. For the convex-hull
support property, see also \citet[Section~4, p.~14]{HanEtAl2026}.

\begin{lemma}[Existence and support of an unrestricted maximizer]
\label{lem:convex-hull-support}
For each sample $x=(x_1,\ldots,x_n)\in(\R^d)^n$, the Gaussian
location-mixture likelihood attains its maximum over
$\mathcal P(\R^d)$. Every maximizer $\widehat G$ satisfies
\[
  \operatorname{supp}(\widehat G)
  \subseteq C=\operatorname{conv}\{x_1,\ldots,x_n\}.
\]
There is a Borel measurable map from $(\R^d)^n$ to $\mathcal P(\R^d)$
that selects a maximizer for each sample, where $\mathcal P(\R^d)$
carries the weak topology.
\end{lemma}

\begin{proof}
Let $\pi_C$ be Euclidean projection onto $C$. Since $x_i\in C$,
\[
  \norm{x_i-\theta}_2^2
  \geq \norm{x_i-\pi_C(\theta)}_2^2
       +\operatorname{dist}(\theta,C)^2.
\]
Thus, for each $i$,
\[
  f_{(\pi_C)_\#G}(x_i)
  =\int_{\R^d}\phi_d(x_i-\pi_C(\theta))\,dG(\theta)
  \geq \int_{\R^d}\phi_d(x_i-\theta)\,dG(\theta)
  =f_G(x_i).
\]
The inequality is strict if $G(C^c)>0$. The log likelihood is continuous on the weakly compact space
$\mathcal P(C)$, so it attains its maximum there. Since projection
onto $C$ never decreases the likelihood, this is also the
unrestricted maximum. Every unrestricted maximizer $\widehat G$
has $\widehat G(C)=1$ by the strict inequality above. Since $C$
is closed, $\operatorname{supp}(\widehat G)\subseteq C$.

For measurable selection, set $r(x)=1+\max_{i,j}|x_{ij}|$ and
$S_x(u)=r(x)u$. Then $S_x$ maps $K_1$ bijectively onto $K_{r(x)}$,
which contains $C$. For each sample $x$, maximize
\[
  \sum_{i=1}^n
  \log\!\left(\int_{K_1}\phi_d(x_i-r(x)u)\,dQ(u)\right)
\]
over $Q$ in the fixed compact metric space $\mathcal P(K_1)$.
The objective is jointly continuous in $(x,Q)$. When the sample
converges, the integrands converge uniformly on $K_1$, and weak
convergence of $Q$ gives convergence of the integrals.

The measurable maximum theorem therefore gives a Borel measurable
choice $x\mapsto Q_x$. Set $\widehat G(x)=(S_x)_\#Q_x$. For every
bounded continuous function $\varphi$,
\[
  \int_{\R^d}\varphi\,d((S_x)_\#Q)
  =\int_{K_1}\varphi(r(x)u)\,dQ(u).
\]
The right-hand side is continuous in $(x,Q)$ by the same
uniform-convergence argument. So $(x,Q)\mapsto(S_x)_\#Q$ is
weakly continuous, and $x\mapsto\widehat G(x)$ is Borel measurable.
It maximizes the likelihood over $\mathcal P(K_{r(x)})$. Since
$C\subseteq K_{r(x)}$, it also attains the unrestricted maximum.
\end{proof}

\begin{proof}[Proof of the upper bound in Theorem~\ref{thm:main}]
First suppose that $G$ is centered and satisfies
\begin{equation}
\label{eq:coordinate-tail-stat}
  G\{|\theta_j|>t\}\leq 2e^{-t^2/(2\sigma^2)},
  \qquad 1\leq j\leq d,\quad t>0.
\end{equation}
Take $A=5$. For $n\geq2$, set
\[
  a_n=\sigma\sqrt{2A\log(2dn)},\qquad
  b_n=\sqrt{2A\log(2dn)},\qquad
  R_n=\max\{1,a_n+b_n\}.
\]
By \eqref{eq:coordinate-tail-stat},
\[
  \tau:=G(K_{a_n}^c)
  \leq\sum_{j=1}^d G\{|\theta_j|>a_n\}
  \leq 2d(2dn)^{-A}<1.
\]
Define $Q=G(\,\cdot\mid K_{a_n})$, and write $p=f_G$ and $q=f_Q$.
Since $\operatorname{TV}(G,Q)=\tau$, contraction under convolution
and the product total variation bound give
\begin{equation}
\label{eq:truncation-tv}
\begin{aligned}
  \operatorname{TV}(p,q)
  &\leq \operatorname{TV}(G,Q)=\tau,\\
  \operatorname{TV}(p^{\otimes n},q^{\otimes n})
  &\leq n\operatorname{TV}(p,q)\leq n\tau,\\
  H^2(p,q)
  &\leq \norm{p-q}_1=2\operatorname{TV}(p,q)\leq2\tau.
\end{aligned}
\end{equation}

Under $q^{\otimes n}$, write $X_i=\Theta_i+Z_i$, where
$\Theta_i\sim Q$ and $Z_i\sim N(0,I_d)$, with all $\Theta_i$ and
$Z_i$ mutually independent. Let
$\mathcal E=\{X_1,\ldots,X_n\in K_{R_n}\}$.
Since $\Theta_i\in K_{a_n}$ almost surely,
\[
  \max_{1\leq i\leq n}\norm{Z_i}_\infty\leq b_n
  \quad\Longrightarrow\quad
  \max_{1\leq i\leq n}\norm{X_i}_\infty
  \leq a_n+b_n\leq R_n.
\]
A union bound over the $dn$ Gaussian coordinates gives
\begin{equation}
\label{eq:sample-compact-event}
  q^{\otimes n}(\mathcal E^c)
  \leq \sum_{i=1}^n\sum_{j=1}^d \Pp\{|Z_{ij}|>b_n\}
  \leq 2dn\,e^{-b_n^2/2}
  =(2dn)^{1-A}.
\end{equation}
Let $\widehat G_n$ be any measurable unrestricted NPMLE and write
$\widehat f_n=f_{\widehat G_n}$. Let $\widehat f_{n,R_n}$ be a
measurable constrained MLE over $\mathcal F_{R_n}$. Such a selection
exists by the measurable maximum theorem. The optimization is over
the fixed compact metric space $\mathcal P(K_{R_n})$, and the log
likelihood is jointly continuous in the sample and the mixing measure.

Since $\mathcal E$ is Borel measurable, the estimator
\[
  \widetilde f_n
  =
  \begin{cases}
    \widehat f_n, & \text{on }\mathcal E,\\
    \widehat f_{n,R_n}, & \text{on }\mathcal E^c
  \end{cases}
\]
is measurable. On $\mathcal E$, the sample convex hull lies in
$K_{R_n}$. Lemma~\ref{lem:convex-hull-support} then gives
$\widehat G_n\in\mathcal P(K_{R_n})$, so $\widehat f_n$ also
maximizes the likelihood over $\mathcal F_{R_n}$. Thus
$\widetilde f_n$ is a constrained MLE for every sample.

By the triangle inequality for $H$ and \eqref{eq:truncation-tv},
\begin{align}
  \E_{p^{\otimes n}}H^2(\widehat f_n,p)
  &\leq 2\E_{p^{\otimes n}}H^2(\widehat f_n,q)+4\tau
  \notag\\
  &\leq 2\E_{q^{\otimes n}}H^2(\widehat f_n,q)
        +4n\tau+4\tau
  \notag\\
  &\leq 2\E_{q^{\otimes n}}H^2(\widetilde f_n,q)
        +4q^{\otimes n}(\mathcal E^c)+4n\tau+4\tau.
  \label{eq:risk-transfer}
\end{align}
The second inequality uses the total variation bound on the
difference of expectations of a loss taking values in $[0,2]$.
The last uses $\widehat f_n=\widetilde f_n$ on $\mathcal E$
and $H^2\leq2$.

Since $Q\in\mathcal P(K_{a_n})\subseteq\mathcal P(K_{R_n})$,
Theorem~\ref{thm:compact} gives
\[
  \E_{q^{\otimes n}}H^2(\widetilde f_n,q)
  \leq C_d\frac{(1+R_n^2+\log n)^d}{n}.
\]
For fixed $d$ and $\sigma$, we have $R_n^2\leq C_{d,\sigma}\log n$
for $n\geq2$. Also, $n\tau\leq(2dn)^{1-A}$ and
$q^{\otimes n}(\mathcal E^c)\leq(2dn)^{1-A}$. Substituting in
\eqref{eq:risk-transfer} and using $H^2\leq2$ yields
\[
  \E_{p^{\otimes n}}H^2(\widehat f_n,p)
  \leq C_{d,\sigma}\min\!\left\{1,\frac{(\log n)^d}{n}\right\}.
\]

For a mixing law $G$ with mean $a$, write $T_b(z)=z+b$ and set
$G_c=(T_{-a})_\#G$. For any mixing laws $M,M'$,
\[
  f_{(T_b)_\#M}(z+b)=f_M(z),\qquad
  H(f_{(T_b)_\#M},f_{(T_b)_\#M'})=H(f_M,f_{M'}).
\]
For centered observations $y_i=X_i-a$, define
\[
  \widehat G_{n,c}(y_1,\ldots,y_n)
  =(T_{-a})_\#\widehat G_n(y_1+a,\ldots,y_n+a).
\]
Translation is continuous on the sample space and on mixing laws in
the weak topology, so $\widehat G_{n,c}$ is measurable. The first
translation identity shows that it maximizes the likelihood for the
centered sample. Since $G_c$ satisfies \eqref{eq:coordinate-tail-stat},
the centered risk bound applies to the right-hand side of
\[
  \E_{f_G^{\otimes n}}H^2(f_{\widehat G_n},f_G)
  =
  \E_{f_{G_c}^{\otimes n}}H^2(f_{\widehat G_{n,c}},f_{G_c}).
\]
This proves \eqref{eq:main}.
\end{proof}

\begin{appendix}
\section{Matching lower bound for each fixed tail parameter}
\label{sec:lower}

The lower-bound construction of \citet[Section~3.3]{KG2022}
perturbs a Gaussian mixing density with covariance
$\min\{\Gamma,1\}I_d$ by products of odd functions, where $\Gamma$
is their tail parameter. When $d\geq2$, these perturbations leave
each coordinate marginal unchanged. Taking $\Gamma=\sigma^2$
therefore gives a matching lower bound for our coordinate tail
class for every fixed $\sigma>0$.

We use the same baseline $N(0,sI_d)$ with $s=\min\{\sigma^2,1\}$
and give a direct proof for both tail classes. Write $\phi_w$ for the one-dimensional
$N(0,w)$ density. Set
\[
  a=\frac{s}{4},
  \qquad
  v^{-1}=\frac{2}{1+a}-\frac{1}{1+s}
        =\frac{1+2s-a}{(1+a)(1+s)}>0.
\]
Define
\[
  q=\frac{av}{1+a+av},
  \qquad
  b=\frac{1+a}{\sqrt{a(1+a+av)}},
  \qquad
  M=\frac{((1+s)v)^{1/4}}{(1+a)^{1/2}}.
\]
Then $0<q<1$. Let $\mathrm{He}_j$ be the probabilists' Hermite
polynomials, defined by
\[
  e^{tz-t^2/2}
  =\sum_{j=0}^{\infty}\mathrm{He}_j(z)\frac{t^j}{j!}.
\]
For $j\geq1$, set
\[
  w_j(u)=M^{-1}q^{-j/2}\phi_a(u)
          \frac{\mathrm{He}_j(bu)}{\sqrt{j!}}.
\]
For $j,k\geq1$, these functions satisfy
\begin{equation}
\label{eq:low-orth}
\begin{aligned}
  \frac{(\phi_1*w_j)(x)}{\sqrt{\phi_{1+s}(x)}}
  &=\sqrt{\phi_v(x)}
    \frac{\mathrm{He}_j(x/\sqrt v)}{\sqrt{j!}},\\
  \left\langle
    \frac{\phi_1*w_j}{\sqrt{\phi_{1+s}}},
    \frac{\phi_1*w_k}{\sqrt{\phi_{1+s}}}
  \right\rangle_{L^2(\R)}
  &=\mathbf{1}_{\{j=k\}}.
\end{aligned}
\end{equation}
Completing the square gives
\[
  \phi_1(x-u)\phi_a(u)
  =\phi_{1+a}(x)\phi_{a/(1+a)}
      \left(u-\frac{ax}{1+a}\right).
\]
For real $t$, Gaussian integration gives
\begin{align*}
  \int_{\R}\phi_1(x-u)\phi_a(u)e^{tbu-t^2/2}\,du
  &=\phi_{1+a}(x)
    \exp\!\left\{\frac{bax}{1+a}t
      -\frac12\left(1-\frac{b^2a}{1+a}\right)t^2\right\}\\
  &=\phi_{1+a}(x)
    \exp\!\left\{\frac{x}{\sqrt v}\sqrt q\,t-\frac12qt^2\right\},
\end{align*}
where $1-b^2a/(1+a)=q$ and $ba/(1+a)=\sqrt{q/v}$.

For $|t|\leq1$, each fixed-order derivative of $e^{tbu-t^2/2}$
is bounded by a polynomial in $|u|$ times $e^{|bu|}$.
The Gaussian factors make this bound integrable. We can therefore
differentiate $j$ times at $t=0$ under the integral. The Hermite
generating function gives
\[
  \int_{\R}\phi_1(x-u)\phi_a(u)\mathrm{He}_j(bu)\,du
  =q^{j/2}\phi_{1+a}(x)\mathrm{He}_j(x/\sqrt v).
\]
Using the definition of $w_j$ and the preceding integral identity,
we obtain
\begin{align*}
  \frac{(\phi_1*w_j)(x)}{\sqrt{\phi_{1+s}(x)}}
  &=\frac{\phi_{1+a}(x)}{M\sqrt{\phi_{1+s}(x)}}
    \frac{\mathrm{He}_j(x/\sqrt v)}{\sqrt{j!}}\\
  &=\sqrt{\phi_v(x)}
    \frac{\mathrm{He}_j(x/\sqrt v)}{\sqrt{j!}}.
\end{align*}
The last equality follows from the definitions of $v$ and $M$.
This proves the first identity in \eqref{eq:low-orth}.

Gaussian integration gives
\[
  \int_{\R}\phi_1(z)e^{tz-t^2/2}e^{rz-r^2/2}\,dz=e^{tr}.
\]
The same Gaussian domination justifies differentiation under the
integral. So
\[
  \int_{\R}\phi_1(z)\mathrm{He}_j(z)\mathrm{He}_k(z)\,dz
  =\left.\partial_t^j\partial_r^k e^{tr}\right|_{(t,r)=(0,0)}
  =j!\,\mathbf{1}_{\{j=k\}}.
\]
By the first identity in \eqref{eq:low-orth} and the substitution
$z=x/\sqrt v$, the inner product in the second identity equals
\begin{align*}
  &\frac{1}{\sqrt{j!\,k!}}
    \int_{\R}\phi_v(x)
      \mathrm{He}_j(x/\sqrt v)\mathrm{He}_k(x/\sqrt v)\,dx\\
  &\quad=\frac{1}{\sqrt{j!\,k!}}
    \int_{\R}\phi_1(z)\mathrm{He}_j(z)\mathrm{He}_k(z)\,dz
    =\mathbf{1}_{\{j=k\}}.
\end{align*}

For $z\in\R$ and $j\geq1$, Cauchy's coefficient estimate on
$|t|=\sqrt j$ gives
\[
  \frac{|\mathrm{He}_j(z)|}{\sqrt{j!}}
  \leq \frac{\sqrt{j!}}{j^{j/2}}e^{|z|\sqrt j+j/2}
  \leq e^{z^2/4+3j/2}.
\]
The last inequality uses $j!\leq j^j$ and
$|z|\sqrt j\leq z^2/4+j$. By the definition of $b$ and
$a=s/4\leq1/4$,
\[
  b^2=\frac{(1+a)^2}{a(1+a+av)}
  <\frac{1+a}{a}
  <\frac{3}{2a}
  =\frac{2}{a}-\frac{2}{s}.
\]
Set $B=M^{-1}\sqrt{s/a}$ and $\Lambda=e^3/q>1$. Then
\begin{align*}
  \frac{|w_j(u)|}{\phi_s(u)}
  &=Bq^{-j/2}
    \exp\!\left\{-\left(\frac{1}{2a}-\frac{1}{2s}\right)u^2\right\}
    \frac{|\mathrm{He}_j(bu)|}{\sqrt{j!}}\\
  &\leq B\Lambda^{j/2}
    \exp\!\left\{-\left(
      \frac{1}{2a}-\frac{1}{2s}-\frac{b^2}{4}
    \right)u^2\right\}\\
  &\leq B\Lambda^{j/2}.
\end{align*}

For sufficiently large $n$, set
\[
  m=\left\lfloor\frac{\log n}{4d\log\Lambda}\right\rfloor\geq1,
  \qquad J=\{1,3,\ldots,2m-1\}^d,
  \qquad \eta_n=\frac{1}{4\sqrt n}.
\]
For $u\in\R^d$ and $\alpha\in\{0,1\}^{J}$, define
\[
  g_0(u)=\prod_{i=1}^d\phi_s(u_i),\qquad
  g_\alpha(u)=g_0(u)+\eta_n\sum_{j\in J}
    \alpha_j\prod_{i=1}^d w_{j_i}(u_i).
\]
Since $|J|=m^d$ and $\sum_{i=1}^d j_i\leq2dm$ for $j\in J$,
\begin{align*}
  \sup_{u\in\R^d}\left|\frac{g_\alpha(u)}{g_0(u)}-1\right|
  &\leq \eta_n \sum_{j\in J}B^d\Lambda^{\frac12\sum_{i=1}^d j_i}\\
  &\leq \eta_n m^dB^d\Lambda^{dm}
   \leq \frac{m^dB^d}{4n^{1/4}}.
\end{align*}
Here $\Lambda^{dm}\leq n^{1/4}$ by the choice of $m$. For fixed
$d$ and $\sigma$, $B$ is constant and $m=O(\log n)$, so the last
bound is at most $1/2$ for all sufficiently large $n$.
The preceding bound gives $\tfrac12 g_0\leq g_\alpha\leq\tfrac32 g_0$.
Each $j_i$ is odd, so $w_{j_i}$ is odd and
$\int_{\R} w_{j_i}(u_i)\,du_i=0$. When $d\geq2$, integrating out
all coordinates except $u_k$ removes every perturbation term. Thus
\[
  \int_{\R^{d-1}}g_\alpha(u)\,du_{-k}
  =\phi_s(u_k),\qquad 1\leq k\leq d.
\]
In particular, $g_\alpha$ is a probability density with mean zero.

Let $G_\alpha$ have density $g_\alpha$ and let $U\sim G_\alpha$.
Since $g_\alpha\leq\tfrac32 g_0$ and $g_0$ is the $N(0,sI_d)$
density, for every unit vector $e$ and every $t>0$,
\[
  \Pp(|e^\top U|>t)
  \leq \tfrac32\Pp(|Z|>t/\sqrt s)
  \leq \tfrac32 e^{-t^2/(2s)}
  \leq 2e^{-t^2/(2\sigma^2)},
  \qquad Z\sim N(0,1).
\]
The Gaussian tail bound follows from
\[
  \Pp(|Z|>z)
  =2e^{-z^2/2}\int_0^\infty\phi_1(y)e^{-zy}\,dy
  \leq e^{-z^2/2},
  \qquad z\geq0.
\]
Thus $G_\alpha$ belongs to the directional-tail class.

Set $f_\alpha=f_{G_\alpha}$ and $f_0(x)=\prod_{i=1}^d\phi_{1+s}(x_i)$.
Convolution gives $\tfrac12 f_0\leq f_\alpha\leq\tfrac32 f_0$ and
\[
  \frac{f_\alpha(x)-f_\beta(x)}{\sqrt{f_0(x)}}
  =\eta_n\sum_{j\in J}(\alpha_j-\beta_j)
    \prod_{i=1}^d
    \frac{(\phi_1*w_{j_i})(x_i)}{\sqrt{\phi_{1+s}(x_i)}}.
\]
The products on the right are orthonormal in $L^2(\R^d)$ by
\eqref{eq:low-orth}. Therefore, with $D(\alpha,\beta)=\#\{j\in J:\alpha_j\ne\beta_j\}$ denoting Hamming distance,
\[
  \int_{\R^d}\frac{(f_\alpha-f_\beta)^2}{f_0}
  =\eta_n^2\sum_{j\in J}(\alpha_j-\beta_j)^2
  =\eta_n^2D(\alpha,\beta).
\]
Since $(\sqrt{f_\alpha}+\sqrt{f_\beta})^2\leq6f_0$,
\[
  H^2(f_\alpha,f_\beta)
  =\int_{\R^d}
    \frac{(f_\alpha-f_\beta)^2}{(\sqrt{f_\alpha}+\sqrt{f_\beta})^2}
  \geq\frac{\eta_n^2}{6}D(\alpha,\beta).
\]
For adjacent vertices $\alpha,\beta$, we have $D(\alpha,\beta)=1$
and $f_\beta\geq f_0/2$. Hence
\[
  \chi^2(f_\alpha\Vert f_\beta)
  =\int_{\R^d}\frac{(f_\alpha-f_\beta)^2}{f_\beta}
  \leq2\int_{\R^d}\frac{(f_\alpha-f_\beta)^2}{f_0}
  =2\eta_n^2.
\]
For adjacent vertices $\alpha$ and $\beta$, Pinsker's inequality,
additivity of relative entropy, and $\mathrm{KL}\leq\chi^2$ give
\begin{align*}
  \operatorname{TV}(f_\alpha^{\otimes n},f_\beta^{\otimes n})
  &\leq \sqrt{\frac n2\,\mathrm{KL}(f_\alpha\Vert f_\beta)}\\
  &\leq \sqrt{\frac n2\,\chi^2(f_\alpha\Vert f_\beta)}
   \leq \sqrt n\,\eta_n=\frac14.
\end{align*}
For any density estimator $\widehat f$, choose $\widehat\alpha$
to minimize $H(\widehat f,f_\beta)$ over $\beta\in\{0,1\}^{J}$,
breaking ties according to a fixed order on $\{0,1\}^{J}$.
By the triangle inequality and the choice of $\widehat\alpha$,
\[
  H(f_{\widehat\alpha},f_\alpha)
  \leq H(f_{\widehat\alpha},\widehat f)+H(\widehat f,f_\alpha)
  \leq 2H(\widehat f,f_\alpha).
\]
The Hellinger separation bound therefore gives
\[
  H^2(\widehat f,f_\alpha)
  \geq \frac14 H^2(f_{\widehat\alpha},f_\alpha)
  \geq \frac{\eta_n^2}{24}D(\widehat\alpha,\alpha).
\]

Write $P_\alpha=f_\alpha^{\otimes n}$. Fix $j\in J$ and pair the
vertices that differ only at $j$. For each pair $(\alpha,\beta)$
with $\alpha_j=0$ and $\beta_j=1$,
\[
  P_\alpha(\widehat\alpha_j=1)+P_\beta(\widehat\alpha_j=0)
  \geq 1-\operatorname{TV}(P_\alpha,P_\beta)
  \geq \frac34.
\]
There are $2^{|J|-1}$ such pairs for each $j$. Summing first over
the pairs and then over $j$ gives
\begin{align*}
  2^{-|J|}\sum_\alpha \E_{P_\alpha}D(\widehat\alpha,\alpha)
  &=\sum_{j\in J}2^{-|J|}
    \sum_\alpha P_\alpha(\widehat\alpha_j\ne\alpha_j)\\
  &\geq \frac{3|J|}{8}.
\end{align*}
The supremum of the risks is at least their average, so
\[
  \sup_\alpha \E_{P_\alpha}H^2(\widehat f,f_\alpha)
  \geq \frac{\eta_n^2}{24}\frac{3|J|}{8}
  =\frac{\eta_n^2|J|}{64}
  =\frac{m^d}{1024n}.
\]
For fixed $d$ and $\sigma$, $m$ is of order $\log n$. Since every
$G_\alpha$ belongs to both tail classes, taking the infimum over
$\widehat f$ proves the lower bound for both classes. When
$\sigma\geq1$, we have $s=1$, so the lower-bound constant depends
only on $d$.
\end{appendix}

\begin{acks}[Acknowledgments]
Xingzhou Wang is also affiliated with the Department of Economics at North Carolina State University.
\end{acks}

\section*{Statement on AI use}
We used ChatGPT (OpenAI) for language editing, reference and formatting
checks, and assistance in checking existing proofs in the article. The authors take responsibility for the
submitted work.

\bibliographystyle{spr-ims-nameyear}
\bibliography{references}
\end{document}